\magnification=1200
%%%%%%%%%%%%%%%%%%%%%%%%%%%%%%%%%%%%%%%%%%%%%%%%%%%%%%%%%%%%%%%%%
%                       \input zaj.sty                          %
\input amssym.def
\input amssym.tex

\hsize=36truecc

\font\secbf=cmb10 scaled 1200
\font\eightrm=cmr8
\font\sixrm=cmr6

\font\eighti=cmmi8

\font\sixi=cmmi6
\skewchar\eighti='177 \skewchar\sixi='177

\font\eightsy=cmsy8
\font\sixsy=cmsy6
\skewchar\eightsy='60 \skewchar\sixsy='60

\font\eightit=cmti8

\font\eightbf=cmbx8
\font\sixbf=cmbx6

\let\sc=\tensc

\font\eightsc=cmcsc10 scaled 800

             % roman do streszcze¤ i literatury
              % italic do literatury

    % Font na tytu'y
 % Font na podtytu'y
\font\secbf=cmb10 scaled 1200
\font\subsecfont=cmb10 scaled \magstephalf
%%%%%%%%%%%%%%%%%%%%%%%%%%%%%%%%%%%%%%%%%%%%%%%%%%
\font\amb=cmmib10

\font\ambi=cmmib10 scaled 700

\newfam\mbfam \def\mb{\textfont1=\amb\fam\mbfam\amb\scriptfont1=\ambi}

\textfont\mbfam\amb \scriptfont\mbfam\ambi

\def\bm#1{\mathchoice
{\hbox{\mb\textfont1=\amb$#1$}}%
{\hbox{\mb\textfont1=\amb$#1$}}%
{\hbox{\mb$\scriptstyle\textfont1=\ambi#1$}}%
{\hbox{\mb$\scriptscriptstyle\textfont1=\ambi#1$}}}
%%%%%%%%%%%%%%%%%%%%%%%%%%%%%%%%%%%%%%%%%%%%%%%%%%%

\def\aa{\def\rm{\fam0\eightrm}%
  \textfont0=\eightrm \scriptfont0=\sixrm \scriptscriptfont0=\fiverm
  \textfont1=\eighti \scriptfont1=\sixi \scriptscriptfont1=\fivei
  \textfont2=\eightsy \scriptfont2=\sixsy \scriptscriptfont2=\fivesy
  \textfont3=\tenex \scriptfont3=\tenex \scriptscriptfont3=\tenex
  \def\sc{\eightsc}
  \def\it{\fam\itfam\eightit}%
  \textfont\itfam=\eightit
  \def\bf{\fam\bffam\eightbf}%
  \textfont\bffam=\eightbf \scriptfont\bffam=\sixbf
   \scriptscriptfont\bffam=\fivebf
  \normalbaselineskip=9.7pt
  \setbox\strutbox=\hbox{\vrule height7pt depth2.6pt width0pt}%
  \normalbaselines\rm}

\def\Proof{\vskip12pt\noindent{\bf Proof.} }

\def\Remark#1{\vskip12pt\noindent{\bf Remark #1}}

\def\m@th{\mathsurround=0pt}

\def\cc#1{\hbox to .89\hsize{$\displaystyle\hfil{#1}\hfil$}\cr}
\def\lc#1{\hbox to .89\hsize{$\displaystyle{#1}\hfill$}\cr}
\def\rc#1{\hbox to .89\hsize{$\displaystyle\hfill{#1}$}\cr}

\def\eqal#1{\null\,\vcenter{\openup\jot\m@th
  \ialign{\strut\hfil$\displaystyle{##}$&&$\displaystyle{{}##}$\hfil
      \crcr#1\crcr}}\,}

\def\section#1{\vskip 22pt plus6pt minus2pt\penalty-400
        {{\secbf
        \noindent#1\rightskip=0pt plus 1fill\par}}
        \par\vskip 12pt plus5pt minus 2pt
        \penalty 1000}

\def\subsection#1{\vskip 20pt plus6pt minus2pt\penalty-400
        {{\subsecfont
        \noindent#1\rightskip=0pt plus 1fill\par}}
        \par\vskip 8pt plus5pt minus 2pt
        \penalty 1000}

\def\subsubsection#1{\vskip 18pt plus6pt minus2pt\penalty-400
        {{\subsecfont
        \noindent#1}}
        \par\vskip 7pt plus5pt minus 2pt
        \penalty 1000}

\def\center#1{{\begingroup \leftskip=0pt plus 1fil\rightskip=\leftskip
\parfillskip=0pt \spaceskip=.3333em \xspaceskip=.5em \pretolerance 9999
\tolerance 9999 \parindent 0pt \hyphenpenalty 9999 \exhyphenpenalty 9999
\par #1\par\endgroup}}

\def\\{\hfill\break}

\def\kwadrat{\hfill$\square$}
\def\mida#1{{{\null\kern-4.2pt\left\bracevert\vbox to 6pt{}\!\hbox{$#1$}\!\right\bracevert\!\!}}}
\def\midy#1{{{\null\kern-4.2pt\left\bracevert\!\!\hbox{$\scriptstyle{#1}$}\!\!\right\bracevert\!\!}}}

\def\diagint{{\raise1.5pt\hbox{$\scriptscriptstyle\diagup$}\hskip-8.7pt\intop}}

\def\divv{{\rm div}\,}
\def\rot{{\rm rot}\,}
\def\const{{\rm const}}

\def\today{${\scriptscriptstyle\number\day-\number\month-\number\year}$}
\footline={{\hfil\rm\the\pageno\hfil${\scriptscriptstyle\rm\jobname}$\ \ \today}}
%%%%%%%%%%%%%%%%%%%%%%%%%%%%%%%%%%%%%%%%%%%%%%%%%%%%%%%%%%%%%%%%%%%%%%%%%%%%%%

\hsize=36truecc
\baselineskip=16truept

\def\D{{\bm{D}}}

\def\I{{\bm{I}}}

\def\T{{\bm{T}}}

\def\0{{\bf0}}

\def\sup{\mathop{\rm sup}\limits}

\def\\{\hfil\break}
\def\N{{\Bbb N}}
\def\R{{\Bbb R}}

\center{\secbf Long time existence of regular solutions to 3d Navier-Stokes 
equations coupled with the heat convection}

\vskip3cm
\centerline{\bf Jolanta Soca\l a$^1$, Wojciech M. Zaj\c aczkowski$^2$}
\vskip1cm
\item{$^1$} State Higher Vocational School in Racib\'orz,\\
S\l owacki Str. 55, 47-400 Racib\'orz, Poland
\item{$^2$}Institute of Mathematics, Polish Academy of Sciences\\
\'Sniadeckich 8, 00-956 Warsaw, Poland\\
and Institute of Mathematics and Cryptology, Cybernetics Faculty,\\
Military University of Technology,\\
Kaliskiego 2, 00-908 Warsaw, Poland\\
e-mail:wz@impan.gov.pl
\vskip1cm

\noindent
{\bf Abstract.} We prove long time existence of regular solutions to the 
Navier-Stokes equations coupled with the heat equation. We consider the system 
in non-axially symmetric cylinder with the slip boundary conditions for the 
Navier-Stokes equations and the Neumann condition for the heat equation. The 
long time existence is possible because we assumed that derivatives with 
respect to the variable along the axis of the cylinder of the initial 
velocity, initial temperature and the external force in $L_2$ norms are 
sufficiently small. We proved the existence of such solutions that velocity 
and temperature belong to $W_\sigma^{2,1}(\Omega\times(0,T))$, where 
$\sigma>{5\over3}$. The existence is proved by the Leray-Schauder fixed point 
theorem.

\noindent
{\bf Key words:} Navier-Stokes equations, heat equation, coupled, slip 
boundary conditions, the Neumann condition, long time existence, regular 
solutions

\noindent
{\bf AMS Subject Classification:} 35D05, 35D10, 35K05, 35K20, 35Q30, 
76D03, 76D05

\vfil\eject

\section{1. Introduction}

We  consider the following problem
$$\eqal{
&v_{,t}+v\cdot\nabla v-\divv\T(v,p)=\alpha(\theta)f\quad &{\rm in}\ \ 
\Omega^T=\Omega\times(0,T),\cr
&\divv v=0\quad &{\rm in}\ \ \Omega^T,\cr
&\theta_{,t}+v\cdot\nabla\theta-\varkappa\Delta\theta=0\quad &{\rm in}\ \ 
\Omega^T,\cr
&\bar n\cdot\D(v)\cdot\bar\tau_\alpha=0,\ \ \alpha=1,2\quad &{\rm on}\ \ 
S^T=S\times(0,T),\cr
&\bar n\cdot\bar v=0\quad &{\rm on}\ \ S^T,\cr
&\bar n\cdot\nabla\theta=0\quad &{\rm on}\ \ S^T,\cr
&v|_{t=0}=v_0,\ \ \theta|_{t=0}=\theta_0\quad &{\rm in}\ \ \Omega,\cr}
\leqno(1.1)
$$
where $\Omega\subset\R^3$ is cylindrical domain, $S=\partial\Omega$, 
$v=(v_1(x,t),v_2(x,t),v_3(x,t))\in\R^3$ is the velocity of the fluid motion, 
$p=p(x,t)\in\R^1$ the pressure, $\theta=\theta(x,t)\in\R_+$ the temperature, 
$f=(f_1(x,t),f_2(x,t),f_3(x,t))\in\R^3$ the external force field, $\bar n$ 
is the unit outward normal vector to the boundary $S$, $\bar\tau_\alpha$, 
$\alpha=1,2$, are tangent vectors to $S$ and the dot denotes the scalar 
product in $\R^3$. We define the stress tensor by
$$
\T(v,p)=\nu\D(v)-p\I,
$$
where $\nu$ is the constant viscosity coefficient, $\I$ is the unit matrix 
and $\D(v)$ is the dilatation tensor of the form
$$
\D(v)=\{v_{i,x_j}+v_{j,x_i}\}_{i,j=1,2,3}.
$$
Finally $\varkappa$ is a positive heat conductivity coefficient.

\noindent
By $x=(x_1,x_2,x_3)$ we denote the Cartesian coordinates, $\Omega\subset\R^3$ 
is a cylindrical type domain parallel to the axis $x_3$ with arbitrary cross 
section.

\noindent
We assume that $S=S_1\cup S_2$, where $S_1$ is the part of the boundary which 
is parallel to the axis $x_3$ and $S_2$ is perpendicular to $x_3$. Hence
$$
S_1=\{x\in\R^3:\ \varphi_0(x_1,x_2)=c_*,\ -b<x_3<b\}
$$
and
$$
S_2=\{x\in\R^3:\ \varphi_0(x_1,x_2)<c_*,\ x_3\ {\rm is\ equal\ either\ to}\ 
-b\ {\rm or}\ b\},
$$
where $b,c_*$ are positive given numbers and $\varphi_0(x_1,x_2)$ describes 
a sufficiently smooth closed curve in the plane $x_3=\const$. We can assume 
$\bar\tau_1=(\tau_{11},\tau_{12},0)$ $\tau_2=(0,0,1)$ and 
$\bar n=(\tau_{12},-\tau_{11},0)$ on $S_1$.

\noindent
Assume $\alpha\in C^2(\R_+)$ and $\Omega^T$ satisfies the weak $l$-horn 
condition, where $l=(2,2,2,1)$.

\noindent
Moreover assume $\Omega$ is not axially symmetric. Now we formulate the main 
result of this paper. Let $g=f_{,x_3}$, $j=v_{,x_3}$, $q=p_{,x_3}$, 
$\vartheta=\theta_{,x_3}$, $\chi=(\rot v)_3$, $F=(\rot f)_3$. Assume that 
$\|\theta(0)\|_{L_\infty(\Omega)}<\infty$.

\noindent
Define
$$
a:[0,\infty)\to[0,\infty),\quad 
a(x)=\sup\{|\alpha(y)|+|\alpha'(y)|+|\alpha''(y)|:|y|\le x\}
$$
and $c_1=a(\|\theta(0)\|_{L_\infty})$. Moreover assume that 
${5\over3}<\sigma<\infty$, ${5\over3}<\varrho<\infty$, 
${5\over\varrho}-{5\over\sigma}<1$ and for $t\le T$
\vskip6pt

\noindent
1. $c_1\|g\|_{L_2(0,t;L_{6/5}(\Omega))}+
c_1c_0\|f\|_{L_\infty(0,t;L_3(\Omega))}+
c_1\|F\|_{L_2(0,t;L_{6/5}(\Omega))}+c_1\|f_3\|_{L_2(0,t;L_{4/3}(S_2))}\break
+\|h(0)\|_{L_2(\Omega)}+\|\vartheta(0)\|_{L_2(\Omega)}+
\|\chi(0)\|_{L_2(\Omega)}+c_0^2(c_1\|f\|_{L_2(0,t;L_{6/5}(\Omega))}+
\|v(0)\|_{L_2(\Omega)})\le k_1<\infty$, 

\noindent
2. $\|f\|_{L_2(0,t;L_3(\Omega))}\le k_2<\infty$,

\noindent
3. $\|f\|_{L_2(\Omega^t)}+\|v_0\|_{H^1(\Omega)}\le k_3<\infty$,

\noindent
4. $c_1\|f\|_{L_\infty(\Omega^t)}e^{cc_1^2k_2^2}k_1+c_1
\|g\|_{L_\sigma(\Omega^t)}+\|\vartheta(0)\|_{W_\sigma^{2-2/\sigma}(\Omega)}+
\|h(0)\|_{W_\sigma^{2-2/\sigma}(\Omega)}\le k_4<\infty$,

\noindent
5. $c_1\|g\|_{L_2(0,t;L_{6/5}(\Omega))}+c_1\|f_3\|_{L_2(0,t;L_{4/3}(S_2))}+
\|h(0)\|_{L_2(\Omega)}+\|\vartheta(0)\|_{L_2(\Omega)}\le d<\infty$,

\noindent
6. $c_1\|f\|_{L_\sigma(\Omega^T)}+\|v(0)\|_{W_\varrho^{2-2/\varrho}(\Omega)}+
\|\theta(0)\|_{W_\varrho^{2-2/\varrho}(\Omega)}\le k_5<\infty$,

\noindent
where $c_0$ is a constant from Lemma 2.2. Assume 
$$
f\in L_\sigma(\Omega^T),\quad g\in L_\sigma(\Omega^T),
$$
$\vartheta(0)\in W_\sigma^{2-2/\sigma}(\Omega)$.

\proclaim Theorem 1.1. 
Let the above assumptions hold. 
Assume that $d$ is sufficiently small (see [5, Main Theorem]).
Then there exists a strong solution $(v,p,\theta)$ to (1.1) such that 
$v,\theta\in W_\varrho^{2,1}(\Omega^T)$, $\nabla p\in L_\varphi(\Omega^T)$, 
$h,\vartheta\in W_\sigma^{2,1}(\Omega^T)$, $\nabla q\in L_\sigma(\Omega^T)$.

The result follows by applying the methods developed in [3] to the more 
complicated system (1.1). However, the proof of existence in this paper 
is much more clear that the one in [3], because the mapping $\phi$ in this 
paper is constructed in a simpler way than the corresponding mapping in [3]. 
This, however, needs more regularity. Therefore in this paper we proved the 
existence of much more regular solutions than in [3].

\section{2. Preliminaries}

In this section we introduce notation and basic estimates for weak solutions 
to problem (1.1).

\subsection{2.1. Notation}

We use isotropic and anisotropic Lebesgue spaces: $L_p(Q)$, 
$Q\in\{\Omega^T,S^T,\Omega,S\}$, $p\in[1,\infty]$; $L_q(0,T;L_p(Q))$ 
$Q\in\{\Omega,S\}$, $p,q\in[1,\infty]$; Sobolev spaces
$$
W_q^{s,s/2}(Q^T),\quad Q\in\{\Omega,S\},\quad q\in[1,\infty],\quad 
s\in\N\cup\{0\}
$$
with the norm
$$
\|u\|_{W_q^{s,s/2}(Q^T)}=\bigg(\sum_{|\alpha|+2a\le s}\intop_{Q^T}
|D_x^\alpha\partial_t^au|^qdxdt\bigg)^{1/q},
$$
where $D_x^\alpha=\partial_{x_1}^{\alpha_1}\partial_{x_2}^{\alpha_2}
\partial_{x_3}^{\alpha_3}$, $|\alpha|=\alpha_1+\alpha_2+\alpha_3$, 
$a,\alpha_i\in\N\cup\{0\}$.

\noindent
In the case $q=2$
$$
H^s(Q)=W_2^s(Q),\quad H^{s,s/2}(Q^T)=W_2^{s,s/2}(Q^T),\quad Q\in\{\Omega,S\}.
$$
Moreover, $L_2(Q)=H^0(Q)$, $L_p(Q)=W_p^0(Q)$, $L_p(Q^T)=W_p^{0,0}(Q^T)$.

\noindent
We define a space natural for study weak solutions to the Navier-Stokes and 
parabolic equations
$$
V_2^k(\Omega^T)=\bigg\{u:\ \|u\|_{V_2^k(\Omega^T)}=
\mathop{\rm esssup}\limits_{t\in[0,T]}
\|u\|_{H^k(\Omega)}+\bigg(\intop_0^T\|\nabla u\|_{H^k(\Omega)}^2dt\bigg)^{1/2}
<\infty\bigg\}.
$$

\subsection{2.2. Weak solutions}

By a weak solution to problem (1.1) we mean $v\in V_2^0(\Omega^T)$, 
$\theta\in V_2^0(\Omega^T)\cap L_\infty(\Omega^T)$ satisfying the integral 
identities
$$\eqal{
&-\intop_{\Omega^T}v\cdot\varphi_{,t}dxdt+\intop_{\Omega^T}v\cdot\nabla v\cdot
\varphi dxdt+{\nu\over2}\intop_{\Omega^T}\D(v)\cdot\D(\varphi)dxdt\cr
&=\intop_{\Omega^T}\alpha(\theta)f\cdot\varphi dxdt+\intop_\Omega v_0\varphi(0)
dx,\cr}
\leqno(2.1)
$$
$$\eqal{
&-\intop_{\Omega^T}\theta\psi_{,t}dxdt+\intop_{\Omega^T}v\cdot\nabla\theta\psi
dxdt+\varkappa\intop_{\Omega^T}\nabla\theta\cdot\nabla\psi dxdt\cr
&=\intop_\Omega\theta_0\psi(0)dx,\cr}
\leqno(2.2)
$$
which hold for $\varphi,\psi\in W_2^{1,1}(\Omega^T)\cap L_5(\Omega^T)$ such 
that $\varphi(T)=0$, $\psi(T)=0$, $\divv\varphi=0$, $\varphi\cdot\bar n|_S=0$.

\proclaim Lemma 2.1. (see [9]) (the Korn inequality) 
Assume that
$$
E_\Omega(v)=\|\D(v)\|_{L_2(\Omega)}^2<\infty,\quad 
v\cdot\bar n|_S=0,\quad \divv v=0.
\leqno(2.3)
$$
If $\Omega$ is not axially symmetric there exists a constant $c_1$ such that
$$
\|v\|_{H^1(\Omega)}^2\le c_1E_\Omega(v).
\leqno(2.4)
$$
If $\Omega$ is axially symmetric, $\eta=(-x_2,x_1,0)$, 
$\alpha=\intop_\Omega v\cdot\eta dx$, then there exists a constant $c_2$ such 
that
$$
\|v\|_{H^1(\Omega)}^2\le c_2(E_\Omega(v)+|\alpha|^2).
\leqno(2.5)
$$

\noindent
Let us consider the problem
$$\eqal{
&h_{,t}-\divv\T(h,q)=f\quad &{\rm in}\ \ \Omega^T,\cr
&\divv h=0\quad &{\rm in}\ \ \Omega^T,\cr
&\bar n\cdot h=0,\ \ \bar n\cdot\D(h)\cdot\bar\tau_\alpha=0,\ \ \alpha=1,2
\quad &{\rm on}\ \ S_1^T,\cr
&h_i=0,\ \ i=1,2,\ \ h_{3,x_3}=0\quad &{\rm on}\ \ S_2^T,\cr
&h|_{t=0}=h(0)\quad &{\rm in}\ \ \Omega\cr}
\leqno(2.6)
$$

\proclaim Theorem 2.1. 
Let $f\in L_p(\Omega^T)$, $h(0)\in W_p^{2-2/p}(\Omega)$, $s\in C^2$, 
$1<p<\infty$. Then there exists a solution to problem (2.6) such that 
$h\in W_p^{2,1}(\Omega^T)$, $\nabla q\in L_p(\Omega^T)$ and there exists 
a constant $c$ depending on $S$ and $p$ such that
$$
\|h\|_{W_p^{2,1}(\Omega^T)}+\|\nabla q\|_{L_p(\Omega^T)}\le c
(\|f\|_{L_p(\Omega^T)}+\|h(0)\|_{W_p^{2-2/p}(\Omega)}).
\leqno(2.7)
$$

\noindent
The proof follows from considerations from [2, Ch. 4].

\noindent
Let us consider the problem
$$\eqal{
&v_{,t}-\divv\T(v,q)=f\quad &{\rm in}\ \ \Omega^T,\cr
&\divv v=0\quad &{\rm in}\ \ \Omega^T,\cr
&\bar n\cdot v=0,\ \ \bar n\cdot\D(v)\cdot\bar\tau_\alpha=0,\ \ \alpha=1,2,
\quad &{\rm on}\ \ S^T,\cr
&v|_{t=0}=v_0\quad &{\rm in}\ \ \Omega.\cr}
\leqno(2.8)
$$

\proclaim Theorem 2.2. (the proof is similar to the proof from [1]) 
Let $f\in L_p(\Omega^T)$, $v(0)\in W_p^{2-2/p}(\Omega)$, $S\in C^2$, 
$1<p<\infty$.
Then there exists a solution to problem (2.8) such that 
$v\in W_p^{2,1}(\Omega^T)$, $\nabla p\in L_p(\Omega^T)$ and there exists 
a constant $c$ depending on $S$ and $p$ such that
$$
\|v\|_{W_p^{2,1}(\Omega^T)}+\|\nabla q\|_{L_p(\Omega^T)}\le c
(\|f\|_{L_p(\Omega^T)}+\|v(0)\|_{W_p^{2-2/p}(\Omega)}).
\leqno(2.9)
$$

\proclaim Lemma 2.3. (see [5, Lemma 2.3]) \\
Assume that $v_0\in L_2(\Omega)$, $\theta_0\in L_\infty(\Omega)$, 
$f\in L_2(0,T;L_{6/5}(\Omega))$, $T<\infty$. Assume that $\Omega$ is not 
axially symmetric. Assume that there exists constants $\theta_*$, $\theta^*$ 
such that $\theta_*<\theta^*$ and
$$
\theta_*\le\theta_0(x)\le\theta^*,\quad x\in\Omega.
$$
Then there exists a weak solution to problem (1.1) such that 
$(v,\theta)\in V_2^0(\Omega^T)\times V_2^0(\Omega^T)$, 
$\theta\in L_\infty(\Omega^T)$ and
$$
\theta_*\le\theta(x,t)\le\theta^*,\quad (x,t)\in\Omega^T,
\leqno(2.10)
$$
$$
\|v\|_{V_2^0(\Omega^T)}\le c(a(\|\theta_0\|_{L_\infty(\Omega)})
\|f\|_{L_2(0,T;L_{6/5}(\Omega))}+\|v_0\|_{L_2(\Omega)})\le c_0,
\leqno(2.11)
$$
$$
\|\theta\|_{V_2^0(\Omega^T)}\le c\|\theta_0\|_{L_2(\Omega)}\le c_0.
\leqno(2.12)
$$

\Remark{2.4.} 
If $\theta(0)\ge0$, then $\theta(t)\ge0$ for $t\ge0$.

\section{3. Existence}

For $\xi,\eta\ge1$, $\sigma,\varrho\ge1$ define
$$\eqal{
&\|(v,\theta)\|_{{\cal M}(\Omega^T)}=\|v\|_{L_\infty(0,T;W_\eta^1(\Omega))}+
\|\theta\|_{L_\infty(0,T;W_\eta^1(\Omega))}\cr
&\quad+\|v_{,x_3}\|L_{\infty(0,T;W_\xi^1(\Omega))}+
\|\theta_{,x_3}\|_{L_\infty(0,T;W_\xi^1(\Omega))},\cr}
$$
$$
{\cal M}(\Omega^T)=\{(v,\theta):\|(v,\theta)\|_{{\cal M}(\Omega^T)}<\infty\},
$$
$$\eqal{
&\|(v,\theta)\|_{{\cal N}(\Omega^T)}=\|v\|_{W_\varrho^{2,1}(\Omega^T)}+
\|\theta\|_{W_\varrho^{2,1}(\Omega^T)}\cr
&\quad+\|v_{,x_3}\|_{W_\sigma^{2,1}(\Omega^T)}+
\|\theta_{,x_3}\|_{W_\sigma^{2,1}(\Omega^T)},\cr}
$$
$$
{\cal N}(\Omega^T)=\{(v,\theta):\|(v,\theta\|_{{\cal N}(\Omega^T)}<\infty\}.
$$

\proclaim Lemma 3.1. 
We have\\
1. $({\cal M}(\Omega^T),\|\ \|_{{\cal M}(\Omega^T)})$ is the Banach space.\\
2. $({\cal N}(\Omega^T),\|\ \|_{{\cal N}(\Omega^T)})$ is the Banach space.\\
3. $\|u\|_{{\cal M}(\Omega^T)}\le c\|u\|_{{\cal N}(\Omega^T)}$ for 
$u\in{\cal N}(\Omega^T)$ and the imbedding 
${\cal N}(\Omega^T)\subset{\cal M}(\Omega^T)$ is compact for $\varrho<\eta$, 
${5\over\varrho}-{3\over\eta}<1$, $\sigma<\xi$, ${5\over\sigma}-{3\over\xi}<1$.

\noindent
Let us consider the problems
$$\eqal{
&v_t-\divv\T(v,p)=-\lambda[\tilde v\cdot\nabla\tilde v+\alpha(\tilde\theta)f],
\cr
&\divv v=0,\cr
&v\cdot\bar n|_S=0,\quad \bar n\cdot\D(v)\cdot\bar\tau_\alpha|_S=0,\quad
\alpha=1,2,\cr
&v|_{t=0}=v_0\cr}
\leqno(3.1)
$$
and
$$\eqal{
&\theta_t-\varkappa\Delta\theta=-\lambda\tilde v\cdot\nabla\tilde\theta,\cr
&\bar n\cdot\nabla\theta|_S=0,\cr
&\theta|_{t=0}=\theta_0,\cr}
\leqno(3.2)
$$
where $\lambda\in[0,1]$ is parameter and $\tilde v,\tilde\theta$ are treated 
as given functions. We will assume that $\alpha\in C^2(\R)$.

\proclaim Lemma 3.2. Assume that
$$\eqal{
&(\tilde v,\tilde\theta)\in{\cal M}(\Omega^T),\quad 3<\eta<\infty,\cr
&f\in L_\varrho(\Omega^T),\quad 1<\varrho<\infty\cr
&v_0\in W_\varrho^{2-2/\varrho}(\Omega),\cr
&S\in C^2,\quad {5\over\varrho}-{3\over\eta}<1,\quad \varrho<\eta.\cr}
$$
Then there exists a unique solution to problem (3.1) such that
$$
v\in W_\varrho^{2,1}(\Omega^T)\subset L_\infty(0,T;W_\eta^1(\Omega))
$$
and
$$\eqal{
&\|v\|_{L_\infty(0,T;W_\eta^1(\Omega))}\le c\|v\|_{W_\varrho^{2,1}(\Omega^T)}
\le c(\lambda\|(\tilde v,\tilde\theta)\|_{{\cal M}(\Omega^T)}^2\cr
&\quad+\lambda a(c\|(\tilde v,\tilde\theta)\|_{{\cal M}(\Omega^T)})
\|f\|_{L_\varrho(\Omega^T)}+\|v_0\|_{W_\varrho^{2-2/\varrho}(\Omega)}).\cr}
$$

\Proof 
We have
$$\eqal{
&\|\tilde v\cdot\nabla\tilde v\|_{L_\varrho(\Omega^T)}\le c
\|\tilde v\|_{L_\infty(\Omega^T)}\|\nabla\tilde v\|_{L_\eta(\Omega^T)}\cr
&\le c\|\tilde v\|_{L_\infty(0,T;W_\eta^1(\Omega))}^2\le c
\|(\tilde v,\tilde\theta)\|_{{\cal M}(\Omega^T)}^2.\cr}
$$
and
$$\eqal{
&\|\alpha(\tilde\theta)f\|_{L_\varrho(\Omega^T)}\le a
(c\|\tilde\theta\|_{L_\infty(0,T;W_\eta^1(\Omega))})
\|f\|_{L_\varrho(\Omega^T)}\cr
&\le c(c\|(\tilde v,\tilde\theta)\|_{{\cal M}(\Omega^T)})
\|f\|_{L_\varrho(\Omega^T)}.\cr}
$$
By Theorem 2.2 the proof is completed.
\kwadrat

\proclaim Lemma 3.3. Assume that
$$\eqal{
&3<\eta<\infty,\quad 1<\varrho<\infty,\quad \varrho<\eta,\quad 
{5\over\varrho}-{3\over\eta}<1,\cr
&(\tilde v,\tilde\theta)\in{\cal M}(\Omega^T),\quad
\theta_0\in W_\varrho^{2-2/\varrho}(\Omega).\cr}
$$
Then there exists a unique solution to problem (3.2) such that
$$
\theta\in W_\varrho^{2,1}(\Omega^T)\subset L_\infty(0,T;W_\eta^1(\Omega))
$$
and
$$%\eqal{&
\|\theta\|_{L_\infty(0,T;W_\eta^1(\Omega))}\le c
\|\theta\|_{W_\varrho^{2,1}(\Omega^T)}%\cr&
\le c(\lambda\|(\tilde v,\tilde\theta)\|_{{\cal M}(\Omega^T)}^2+
\|\theta_0\|_{W_\varrho^{2-2/\varrho}(\Omega)}).%\cr}
$$

\Proof 
We have
$$\eqal{
&\|\tilde v\cdot\nabla\tilde\theta\|_{L_\varrho(\Omega^T)}\le
\|\tilde v\|_{L_\infty(\Omega^T)}\|\nabla\tilde\theta\|_{L_\eta(\Omega^T)}\cr
&\le c\|\tilde v\|_{L_\infty(0,T;W_\eta^1(\Omega))}
\|\tilde\theta\|_{L_\infty(0,T;W_\eta^1(\Omega))}\cr
&\le c\|(\tilde v,\tilde\theta)\|_{{\cal M}(\Omega^T)}^2.\cr}
$$
Then similarly as in Theorem 9.1 from [2, Ch. 4, Sect. 9] (see also 
[6, Theorem 17]) we prove the lemma.
\kwadrat

\proclaim Lemma 3.4. Let
$$\eqal{
&(\tilde v,\tilde\theta)\in{\cal M}(\Omega^T),\quad 3<\xi<\infty,\quad 
3<\eta<\infty\cr
&f\in L_\sigma(\Omega^T),\quad g\in L_\sigma(\Omega^T),\quad 
1<\sigma<\infty\quad (where\ \ g=f_{,x_3})\cr
&\sigma<\eta,\quad S\in C^2,\quad \sigma<\xi,\quad 
{5\over\sigma}-{3\over\xi}<1.\cr}
$$
Let $v,p$ be a unique solution to problem (3.1).
Let $h=v_{,x_3}$, $q=p,x_3$. Assume $h(0)\in W_\sigma^{2-2/\sigma}(\Omega)$. 
Then
$$
h\in W_\sigma^{2,1}(\Omega^T)\subset L_\infty(0,T;W_\xi^1(\Omega))
$$
and
$$\eqal{
&\|h\|_{L_\infty(0,T;W_\xi^1(\Omega))}\le c\|h\|_{W_\sigma^{2,1}(\Omega^T)}
\le c(\lambda\|(\tilde v,\tilde\theta)\|_{{\cal M}(\Omega)}^2\cr
&\quad+\lambda a(c\|(\tilde v,\tilde\theta)\|_{{\cal M}(\Omega^T)})
\|(\tilde v,\tilde\theta)\|_{{\cal M}(\Omega^T)}\|f\|_{L_\sigma(\Omega^T)}\cr
&\quad+\lambda a(c\|(\tilde v,\tilde\theta)\|_{{\cal M}(\Omega^T)})
\|g\|_{L_\sigma(\Omega^T)}+\|h(0)\|_{W_\sigma^{2-2/\sigma}(\Omega)}).\cr}
$$

\Proof 
The function $h$ is a solution of the following problem
$$\eqal{
&h_{,t}-\divv\T(h,q)=\lambda[-\tilde v\cdot\nabla\tilde h-\tilde h\cdot\nabla
\tilde v+\alpha_\theta(\tilde\theta)\tilde vf+\alpha(\tilde\theta)g]\quad 
&{\rm in}\ \ \Omega^T,\cr
&\divv h=0\quad &{\rm in}\ \ \Omega^T,\cr
&\bar n\cdot h=0,\ \ \bar n\cdot\D(h)\cdot\bar\tau_\alpha=0,\ \ \alpha=1,2\quad
&{\rm on}\ \ S_1^T,\cr
&h_i=0,\ \ i=1,2,\ \ h_{3,x_3}=0\quad &{\rm on}\ \ S_2^T,\cr
&h|_{t=0}=h(0)\quad &{\rm in}\ \ \Omega,\cr}
$$
where $\tilde h=\tilde v_{,x_3}$, $\tilde\vartheta=\tilde\theta_{,x_3}$.

\noindent
We have
$$\eqal{
&\|\tilde v\cdot\nabla\tilde h\|_{L_\sigma(\Omega^T)}\le c
\|\tilde v\|_{L_\infty(\Omega^T)}\|\nabla\tilde h\|_{L_\xi(\Omega^T)}\cr
&\le c\|\tilde v\|_{L_\infty(0,T;W_\eta^1(\Omega))}
\|\tilde h\|_{L_\infty(0,T;W_\xi^1(\Omega))}\cr
&\le c\|(\tilde v,\tilde\theta)\|_{{\cal M}(\Omega^T)}^2.\cr}
$$
and
$$\eqal{
&\|\tilde h\cdot\nabla\tilde v\|_{L_\sigma(\Omega^T)}\le c
\|\tilde h\|_{L_\infty(\Omega^T)}\|\nabla\tilde v\|_{L_\eta(\Omega^T)}\cr
&\le c\|\tilde h\|_{L_\infty(0,T;W_\xi^1(\Omega))}
\|\tilde v\|_{L_\infty(0,T;W_\eta^1(\Omega))}\cr
&\le c\|(\tilde v,\tilde\theta)\|_{{\cal M}(\Omega^T)}^2.\cr}
$$
Then
$$\eqal{
&\|\alpha_\theta(\tilde\theta)\tilde\vartheta f\|_{L_\sigma(\Omega^T)}\le ca
(c\|\tilde\theta\|_{L_\infty(0,T;W_\eta^1(\Omega))})
\|\tilde\vartheta\|_{L_\infty(0,T;W_\xi^1(\Omega))}
\|f\|_{L_\sigma(\Omega^T)}\cr
&\le ca(c\|(\tilde v,\tilde\theta)\|_{{\cal M}(\Omega^T)})
\|(\tilde v,\tilde\theta)\|_{{\cal M}(\Omega^)}\|f\|_{L_\sigma(\Omega^T)}.\cr}
$$
We have
$$\eqal{
&\|\alpha(\tilde\theta)g\|_{L_\sigma(\Omega^T)}\le a(c
\|\tilde\theta\|_{L_\infty(0,T;W_\eta^1(\Omega))})\|g\|_{L_\sigma(\Omega^T)}
\cr
&\le c(c\|(\tilde v,\tilde\theta)\|_{{\cal M}(\Omega^T)})
\|g\|_{L_\sigma(\Omega^T)}.\cr}
$$
By Theorem 2.1 the proof is completed.
\kwadrat

\proclaim Lemma 3.5. 
Assume that $3<\eta<\infty$, $1<\sigma<\infty$, $\sigma<\eta$, 
${5\over\sigma}-{3\over\xi}<1$, $3<\xi<\infty$, $\sigma<\xi$, 
$(\tilde v,\tilde\theta)\in{\cal M}(\Omega^T)$.
Let $\theta$ be a unique solution to problem (3.2). 
Let $\vartheta=\theta_{,x_3}$. 
Assume that $\vartheta(0)\in W_\sigma^{2-2/\sigma}(\Omega)$. Then
$$
\vartheta\in W_\sigma^{2,1}(\Omega^T)\subset L_\infty(0,T;W_\xi^1(\Omega))
$$
and
$$
\|\vartheta\|_{L_\infty(0,T;W_\xi^1(\Omega))}\le c
\|\vartheta\|_{W_\sigma^{2,1}(\Omega^T)}\le c(\lambda
\|(\tilde v,\tilde\theta)\|_{{\cal M}(\Omega^T)}^2
+\|\vartheta(0)\|_{W_\sigma^{2-2/\sigma}(\Omega)}).
$$

\Proof 
The function $\vartheta$ is solution of the following problem
$$\eqal{
&\vartheta_{,t}-\varkappa\Delta\vartheta=-\lambda[\tilde h\cdot\nabla
\tilde\theta+\tilde v\cdot\nabla\tilde\vartheta]\quad &{\rm in}\ \ \Omega^T,\cr
&\bar n\cdot\nabla\vartheta=0\quad &{\rm on}\ \ S_1^T,\cr
&\vartheta=0\quad &{\rm on}\ \ S_2^T,\cr
&\vartheta|_{t=0}=\vartheta(0)\quad &{\rm in}\ \ \Omega,\cr}
$$
where $\tilde\vartheta=\tilde\theta_{,x_3}$. We have
$$\eqal{
&\|\tilde h\cdot\nabla\tilde\theta\|_{L_\sigma(\Omega^T)}\le
\|\tilde h\|_{L_\infty(\Omega^T)}\|\nabla\tilde\theta\|_{L_\eta(\Omega^T)}\cr
&\le c\|\tilde h\|_{L_\infty(0,T;W_\xi^1(\Omega))}
\|\tilde\theta\|_{L_\infty(0,T;W_\eta^1(\Omega))}\cr
&\le c\|(\tilde v,\tilde\theta)\|_{{\cal M}(\Omega^T)}^2\cr}
$$
and
$$\eqal{
&\|\tilde v\cdot\nabla\tilde\vartheta\|_{L_\sigma(\Omega^T)}\le
\|\tilde v\|_{L_\infty(\Omega^T)}\|\nabla\tilde\vartheta\|_{L_\xi(\Omega^T)}\cr
&\le c\|\tilde v\|_{L_\infty(0,T;W_\eta^1(\Omega))}
\|\tilde\vartheta\|_{L_\infty(0,T;W_\xi^1(\Omega))}\cr
&\le c\|(\tilde v,\tilde\theta)\|_{{\cal M}(\Omega^T)}^2.\cr}
$$
Then similarly as in Theorem 9.1 from [2, Ch. 4, Sect. 9] (see also 
[6, Theorem 17]) we prove the lemma.
\kwadrat

\noindent
From Lemmas 3.1--3.5 it follows that if 
$(\tilde v,\tilde\theta)\in{\cal M}(\Omega^T)$, then there exists a unique 
solution $(v,\theta)$ to problems (3.1)--(3.2) such that 
$(v,\theta)\in{\cal M}(\Omega^T)$.

To prove the existence of solutions to problem (1.1) we apply the 
Leray-Schauder fixed point theorem (see [4, 7, 8]). Therefore we introduce the 
mapping $\phi:[0,1]\times{\cal M}(\Omega^T)\to{\cal M}(\Omega^T)$, 
$(\lambda,\tilde v,\tilde\theta)\to\phi(\lambda,\tilde v,\tilde\theta)=
(v,\theta)$ where $(v,\theta)$ is a solution to problems (3.1)--(3.2).

\noindent
For $\lambda=0$ we have the existence of a unique solution. For $\lambda=1$ 
every fixed point is a solution to problem (1.1).

\proclaim Lemma 3.6. 
Let the assumptions of Lemmas 3.2--3.5 be satisfied. Then the mappings 
$\phi(\lambda,\cdot):{\cal M}(\Omega^T)\to{\cal M}(\Omega^T)$, 
$\lambda\in[0,1]$ are completely continuous.

\Proof 
By Lemmas 3.1--3.5 the mappings $\phi(\lambda,\cdot)$, $\lambda\in[0,1]$ are 
compact. From this it follows that bounded sets in ${\cal M}(\Omega^T)$ 
are transformed into bounded sets in ${\cal M}(\Omega^T)$. Let 
$(\tilde v_i,\tilde\theta_i)\in{\cal M}(\Omega^T)$, $i=1,2$ be two given 
elements. Then $(v_i,\theta_i)$, $i=1,2$ are solutions to the problems
$$\eqal{
&v_{it}-\divv\T(v_i,p_i)=
-\lambda(\tilde v_i\cdot\nabla\tilde v_i+\alpha(\tilde\theta_i)f),\cr
&\divv v_i=0\cr
&\bar n\cdot\D(v_i)\cdot\bar\tau|_S=0,\ \ \bar n\cdot v_i|_S=0,\cr
&v_i|_{t=0}=v_0,\ \ i=1,2,\cr}
\leqno(3.3)
$$
and
$$\eqal{
&\theta_{it}-\varkappa\Delta\theta_i=-\lambda\tilde v_i\cdot\nabla
\tilde\theta_i,\cr
&\bar n\cdot\nabla\theta_i|_S=0,\cr
&\theta_i|_{t=0}=\theta_0,\ \ i=1,2.\cr}
\leqno(3.4)
$$
To show continuity we introduce the differences
$$
V=v_1-v_2,\quad P=p_1-p_2,\quad {\cal T}=\theta_1-\theta_2
\leqno(3.5)
$$
which are solutions to the problems
$$\eqal{
&V_t-\divv\T(V,P)=-\lambda[\tilde V\cdot\nabla\tilde v_1+\tilde v_2\cdot
\nabla\tilde V+(\alpha(\tilde\theta_1)-\alpha(\tilde\theta_2))f]\cr
&\divv V=0\cr
&V\cdot\bar n|_S=0\ \ \bar n\cdot\D(V)\cdot\bar\tau|_S=0,\cr
&V|_{t=0}=0\cr}
\leqno(3.6)
$$
and
$$\eqal{
&{\cal T}_t-\varkappa\Delta{\cal T}=-\lambda[\tilde V\cdot\nabla\tilde\theta_1+
\tilde v_2\cdot\nabla\tilde{\cal T}]\cr
&\bar n\cdot\nabla{\cal T}|_S=0,\cr
&{\cal T}|_{t=0}=0,\cr}
\leqno(3.7)
$$
where $\tilde V=\tilde v_1-\tilde v_2$, 
$\tilde{\cal T}=\tilde\theta_1-\tilde\theta_2$.

\noindent
In view of [3] and [7, 8] we have
$$\eqal{
&\|V\|_{W_\varrho^{2,1}(\Omega^T)}+\|{\cal T}\|_{W_\varrho^{2,1}(\Omega^T)}
\le c
[\|\tilde V\|_{L_\infty(\Omega^T)}\|\nabla\tilde v_1\|_{L_\varrho(\Omega^T)}\cr
&\quad+\|\tilde v_2\|_{L_\infty(\Omega^T)}
\|\nabla\tilde V\|_{L_\varrho(\Omega^T)}+ca(\max\{
\|\tilde\theta_1\|_{L_\infty(\Omega^T)},
\|\tilde\theta_2\|_{L_\infty(\Omega^T)}\})\cdot\cr
&\quad\cdot\|\tilde{\cal T}\|_{L_\infty(\Omega^T)}\|f\|_{L_\varrho(\Omega^T)}+
\|\tilde v_2\|_{L_\infty(\Omega^T)}
\|\nabla\tilde{\cal T}\|_{L_\varrho(\Omega^T)}\cr
&\quad+\|\tilde V\|_{L_\infty(\Omega^T)}
\|\nabla\tilde\theta_1\|_{L_\varrho(\Omega^T)}]\le c
(\|\tilde V\|_{{\cal M}(\Omega^T)}+\|\tilde{\cal T}\|_{{\cal M}(\Omega^T)}).
\cr}
\leqno(3.8)
$$
Let $h_i=v_{i,x_3}$, $q_i=p_{i,x_3}$, $\vartheta_i=\theta_{i,x_3}$, 
$\tilde h_i=\tilde v_{i,x_3}$, $\tilde\vartheta_i=\tilde\theta_{i,x_3}$.

\noindent
The functions $h_i,\vartheta_i$, $i=1,2$ are solutions to the following 
problems
$$\eqal{
&h_{i,t}-\divv\T(h_i,q_i)=-\lambda[\tilde h_i\cdot\nabla\tilde v_i+
\tilde v_i\cdot\nabla\tilde h_i+\alpha_\theta(\tilde\theta_i)
\tilde\vartheta_if+\alpha(\tilde\theta_i)g]\quad &{\rm in}\ \ \Omega^T,\cr
&\divv h_i=0\quad &{\rm in}\ \ \Omega^T,\cr
&\bar n\cdot h_i=0,\ \ \bar n\cdot\D(h)\cdot\bar\tau_\alpha,\ \ \alpha=1,2,\ \ 
i=1,2\quad &{\rm on}\ \ S_1^T,\cr
&h_{ij}=0,\ \ i=1,2,\ \ j=1,2\quad &{\rm on}\ \ S_2^T,\cr
&h_{i3,x_3}=0,\ \ i=1,2\quad &{\rm on}\ \ S_2^T,\cr
&h_i|_{t=0}=h(0)\quad &{\rm in}\ \ \Omega\cr}
$$
and
$$\eqal{
&\vartheta_{i,t}-\varkappa\Delta\vartheta_i=-\lambda[\tilde h_i\cdot\nabla
\tilde\theta_i+\tilde v_i\cdot\nabla\tilde\vartheta_i]\quad &{\rm in}\ \ 
\Omega^T,\cr
&\bar n\cdot\nabla\vartheta_i=0\quad &{\rm on}\ \ S_1^T,\cr
&\vartheta_i=0\quad &{\rm on}\ \ S_2^T,\cr
&\vartheta_i|_{t=0}=\vartheta(0)\quad &{\rm in}\ \ \Omega.\cr}
$$
We introduce the differences
$$
H=h_1-h_2,\quad Q=q_1-q_2,\quad R=\vartheta_1-\vartheta_2
$$
which are solutions to the problems
$$\eqal{
&H_{,t}-\divv\T(H,Q)=-\lambda[\tilde H\cdot\nabla\tilde v_1+\tilde h_2\cdot
\nabla\tilde V+\tilde V\cdot\nabla\tilde h_1+\tilde v_2\cdot\nabla\tilde H
\quad\cr
&\quad+(\alpha_\theta(\tilde\theta_1)-\alpha_\theta(\tilde\theta_2))
\tilde\vartheta_1f+\alpha_\theta(\tilde\theta_2)\tilde Rf\cr
&\quad+(\alpha(\tilde\theta_1)-\alpha(\tilde\theta_2))g\quad &{\rm in}\ \ 
\Omega^T,\cr
&\divv H=0\quad &{\rm in}\ \ \Omega^T,\cr
&\bar n\cdot H=0,\ \ \bar n\cdot\D(H)\cdot\bar\tau_\alpha=0,\ \ \alpha=1,2\quad
&{\rm on}\ \ S_1^T,\cr
&H_j=0,\ \ j=1,2,\ \ H_{3,x_3}=0\quad &{\rm on}\ \ S_2^T,\cr
&H|_{t=0}=0\quad &{\rm in}\ \ \Omega.\cr}
$$
and
$$\eqal{
&R_{,t}-\varkappa\Delta R=-\lambda[\tilde H\cdot\nabla\tilde\theta_1+
\tilde h_2\cdot\nabla\tilde{\cal T}+\tilde V\cdot\nabla\tilde\vartheta_1
+\tilde v_2\cdot\nabla\tilde R]\quad &{\rm in}\ \ \Omega^T,\cr
&\bar n\cdot\nabla R=0\quad &{\rm on}\ \ S_1^T,\cr
&R=0\quad &{\rm on}\ \ S_2^T,\cr
&R|_{t=0}=0\quad &{\rm in}\ \ \Omega,\cr}
$$
where $\tilde H=\tilde h_1-\tilde h_2$, 
$\tilde R=\tilde\vartheta_1-\tilde\vartheta_2$. In view of [3] and [7, 8] 
we have
$$\eqal{
&\|H\|_{W_\sigma^{2,1}(\Omega^T)}+\|R\|_{W_\sigma^{2,1}(\Omega^T)}\le c
[\|\tilde H\|_{L_\infty(\Omega^T)}\|\nabla\tilde v_1\|_{L_\eta(\Omega^T)}\cr
&\quad+\|\tilde h_2\|_{L_\infty(\Omega^T)}\|\nabla\tilde V\|_{L_\eta(\Omega^T)}
+\|\tilde V\|_{L_\infty(\Omega^T)}\|\nabla\tilde h_1\|_{L_\xi(\Omega^T)}\cr
&\quad+\|\tilde v_2\|_{L_\infty(\Omega^T)}
\|\nabla\tilde H\|_{L_\xi(\Omega^T)}+c(\max
\{\|\tilde\theta_1\|_{L_\infty(\Omega^T)},
\|\tilde\theta_2\|_{L_\infty(\Omega^T)}\})\cr
&\quad\cdot\|{\cal T}\|_{L_\infty(\Omega^T)}
\|\tilde\vartheta_1\|_{L_\infty(\Omega^T)}\|f\|_{L_\sigma(\Omega^T)}\cr
&\quad+c(\|\tilde\theta_2\|_{L_\infty(\Omega^T)}
\|\tilde R\|_{L_\infty(\Omega^T)}\|f\|_{L_\sigma(\Omega^T)}\cr
&\quad+a(\max\{\|\tilde\theta_1\|_{L_\infty(\Omega^T)},
\|\tilde\theta_2\|_{L_\infty(\Omega^T)}\})
\|\tilde{\cal T}\|_{L_\infty(\Omega^T)}\|g\|_{L_\sigma(\Omega^T)}\cr
&\quad+\|\tilde H\|_{L_\infty(\Omega^T)}
\|\nabla\tilde\theta_1\|_{L_\eta(\Omega^T)}+\|\tilde h_2\|_{L_\infty(\Omega^T)}
|\nabla\tilde{\cal T}\|_{L_\eta(\Omega^T)}\cr
&\quad+\|\tilde V\|_{L_\infty(\Omega^T)}
\|\nabla\tilde\vartheta_1\|_{L_\xi(\Omega^T)}+
\|\tilde v_2\|_{L_\infty(\Omega^T)}\|\nabla\tilde R\|_{L_\xi(\Omega^T)}
\le c(\|\tilde V,\tilde{\cal T}\|_{{\cal M}(\Omega^T)})\cr}
$$
and from (3.8) and Lemma 3.1 we obtain
$$
\|(V,{\cal T})\|_{{\cal M}(\Omega^T)}\le c
\|(\tilde V,\tilde{\cal T})\|_{{\cal M}(\Omega^T)}.
$$
Hence continuity of $\phi$ follows. This concludes the proof.
\kwadrat

\proclaim Lemma 3.7. 
Let assumptions of Lemmas 3.2--3.5 be satisfied. Then for every bounded 
subset ${\cal M}_0$ of ${\cal M}(\Omega^T)$, the family of maps
$$
\phi(\cdot,\tilde v,\tilde\theta):[0,1]\to{\cal M}(\Omega^T),\quad 
(\tilde v,\tilde\theta)\in{\cal M}_0
$$
is uniformly equicontinuous.

\Proof 
Let $(\tilde v,\tilde\theta)\in{\cal M}_0$, $\lambda_i\in[0,1]$, $i=1,2$, 
$\lambda_1\ge\lambda_2$ and $v_i,\theta_i$ are solutions to the problems
$$\eqal{
&v_{it}-\divv\T(v_i,p_i)=-\lambda_i(\tilde v\cdot\nabla\tilde v+
\alpha(\tilde\theta)f),\cr
&\divv v_i=0,\cr
&\bar n\cdot\D(v_i)\cdot\bar\tau|_S=0,\ \ \bar n\cdot v_i|_S=0,\cr
&v_i|_{t=0}=v_0,\ \ i=1,2\cr}
$$
and
$$\eqal{
&\theta_{it}-\varkappa\Delta\theta_i=-\lambda_i\tilde v\cdot\nabla\tilde\theta,
\cr
&\bar n\cdot\nabla\theta_i|_S=0,\cr
&\theta_i|_{t=0}=\theta_0,\ \ i=1,2.\cr}
$$
To show uniform equicontinuity we introduce the differences
$$
V=v_1-v_2,\quad P=p_1-p_2,\quad {\cal T}=\theta_1-\theta_2
$$
which are solutions to the problems
$$\eqal{
&V_t-\divv\T(V,P)=-(\lambda_1-\lambda_2)(\tilde v\cdot\nabla\tilde v+\alpha
(\tilde\theta)f),\cr
&\divv V=0,\cr
&\bar n\cdot\D(V)\cdot\bar\tau|_S=0,\ \ \bar n\cdot V|_S=0,\cr
&V|_{t=0}=0\cr}
$$
and
$$\eqal{
&{\cal T}_{,t}-\varkappa\Delta{\cal T}=-(\lambda_1-\lambda_2)\tilde v\cdot
\nabla\tilde\theta,\cr
&\bar n\cdot\nabla{\cal T}|_S=0,\cr
&{\cal T}|_{t=0}=0.\cr}
$$
In view of Lemmas 3.2--3.3
$$\eqal{
&\|V\|_{L_\infty(0,T;W_\eta^1(\Omega))}+
\|{\cal T}\|_{L_\infty(0,T;W_\eta^1(\Omega))}\le c((\lambda_1-\lambda_2)
\|(\tilde v,\tilde\theta)\|_{{\cal M}(\Omega)}^2\cr
&\quad+(\lambda_1-\lambda_2)a(c
\|(\tilde v,\tilde\theta)\|_{{\cal M}(\Omega^T)})\|f\|_{L_\varphi(\Omega^T)}.
\cr}
\leqno(3.9)
$$
Let $h_i=v_{i,x_3}$, $\vartheta_i=\theta_{i,x_3}$.

\noindent
We introduce the differences
$$
H=h_1-h_2,\quad R=\vartheta_1-\vartheta_2
$$
which satisfy the following conditions
$$
H=V_{,x_3},\quad R={\cal T}_{,x_3}.
$$
In view of Lemmas 3.4 and 3.5
$$\eqal{
&\|H\|_{L_\infty(0,T;W_\xi^1(\Omega))}+\|R\|_{L_\infty(0,T;W_\xi^1(\Omega))}
\le c((\lambda_1-\lambda_2)\|(\tilde v,\tilde\theta)\|_{{\cal M}(\Omega^T)}\cr
&\quad+(\lambda_1-\lambda_2)a
(c\|(\tilde v,\tilde\theta)\|_{{\cal M}(\Omega^T)})
\|(\tilde v,\tilde\theta)\|_{{\cal M}(\Omega^T)}\|f\|_{L_\sigma(\Omega^T)}\cr
&\quad+(\lambda_1-\lambda_2)a(c
\|(\tilde v,\theta)\|_{{\cal M}(\Omega^T)}\|g\|_{L_\sigma(\Omega^T)}.\cr}
\leqno(3.10)
$$
From (3.9) and (3.10) the uniform equicontinuity of 
$\phi(\cdot,\tilde v,\tilde\theta)$ follows.

\noindent
{\bf Proof of Theorem 1.1.}\\
In view of the above considerations and [5, Main Theorem] the assumptions of 
the Leray-Schauder fixed point theorem are satisfied. Hence Theorem 1.1 is 
proved.

\section{References}

\item{1.} Alame, W.: On existence of solutions for the nonstationary 
Stokes system with slip boundary conditions, Appl. Math. 32 (2) (2005), 
195--223.

\item{2.} Ladyzhenskaya, O. A.; Solonnikov, V. A.; Uraltseva, N. N.: Linear 
and quasilinear equations of parabolic type, Nauka, Moscow 1967 (in Russian).

\item{3.} Renc\l awowicz, J.; Zaj\c aczkowski, W. M.: Large time regular 
solutions to the Navier\--Stokes equations in cylindrical domains, TMNA 32 
(2008), 69--87.

\item{4.} Soca\l a, J.; Zaj\c aczkowski, W. M.: Long time existence of 
solutions to 2d Navier-Stokes equations with heat convection, Appl. Math. 36 
(4) (2009), 453--463.

\item{5.} Soca\l a, J.; Zaj\c aczkowski, W. M.: Long time estimate of 
solutions to 3d Navier-Stokes equations coupled with the heat convection.

\item{6.} Solonnikov, V. A.: A priori estimates for second order parabolic 
equations, Trudy MIAN 70 (1964), 133--212 (in Russian).

\item{7.} Zaj\c aczkowski, W. M.: Long time existence of regular solutions to 
Navier-Stokes equations in cylindrical domains under boundary slip conditions, 
Studia Math. 169 (3) (2005), 243--285.

\item{8.} Zaj\c aczkowski, W. M.: Global special solutions to the 
Navier-Stokes equaitons in a cylindrical domain without the axis of symmetry, 
TMNA 24 (2004), 69--105.

\item{9.} Zaj\c aczkowski, W. M.: Global existence of axially symmetric 
solutions of incompressible Navier-Stokes equations with large angular 
component of velocity, Colloq, Math. 100 (2004), 243--263.

\bye